\documentclass[reqno,12pt, oneside]{amsart}

\usepackage{enumerate}
\usepackage[nobysame]{amsrefs}
\usepackage{enumitem}
\usepackage{amssymb}
\usepackage{amsmath}
\usepackage{amsthm}
\usepackage{mathtools}
\usepackage{amscd}
\usepackage{microtype}
\usepackage[right=2cm,left=2cm, top=2.5cm, bottom=2.5cm,includehead]{geometry}
\usepackage{amsfonts}

\usepackage{tgpagella}
\usepackage{eulervm}

\theoremstyle{plain}
\newtheorem{theorem}{Theorem}

\newtheorem{proposition}[theorem]{Proposition}

\theoremstyle{remark}
\newtheorem{remark}[theorem]{Remark}

\theoremstyle{definition}

\DeclarePairedDelimiter{\abs}{\lvert}{\rvert}

\DeclarePairedDelimiter{\norm}{\lVert}{\rVert}

\DeclarePairedDelimiterX{\dset}[2]{\lbrace}{\rbrace}{#1\;\delimsize|\;#2}

\DeclareMathOperator{\var}{var}

\newcommand{\ball}{{\overline{B}}}

\newcommand{\f}{\varphi}

\title[Integral operators in Schramm bounded variation spaces]{Integral operators in the spaces of functions of bounded Schramm variation}

\author{Jacek Gulgowski}
\address[J.~Gulgowski]{Institute of Mathematics\\
Faculty of Mathematics, Physics and Informatics\\
University of Gda\'nsk\\
80-308 Gda\'nsk\\
Poland}
\email[J.~Gulgowski]{dzak@mat.ug.edu.pl}

\author{Piotr Kasprzak}
\address[P. Kasprzak]{Department of Nonlinear Analysis and Applied Topology\\
Faculty of Mathematics and Computer Science\\
  Adam Mickiewicz University in Pozna\'n\\
  ul.\ Uniwersytetu Pozna\'nskiego 4\\
  61-614 Pozna\'n\\
  Poland}
\email[P.~Kasprzak]{kasp@amu.edu.pl}

\author{Piotr Ma\'ckowiak}
\address[P. Ma\'ckowiak]{Department of Nonlinear Analysis and Applied Topology\\
Faculty of Mathematics and Computer Science\\
  Adam Mickiewicz University in Pozna\'n\\
  ul.\ Uniwersytetu Pozna\'nskiego 4\\
  61-614 Pozna\'n\\
  Poland}
\email[P.~Ma\'ckowiak]{piotr.mackowiak@amu.edu.pl}

\keywords{compact linear operator, compactness criterion, equinormed set, functions of bounded Schramm variation, linear integral operators, precompact set, relatively compact set}
\subjclass[2010]{47G10, 46B50, 26A45}
\date{\today}

\begin{document}
\begin{abstract}
In this paper we provide a full characterization of linear integral operators acting from the space of functions of bounded Jordan variation to the space of functions of bounded Schramm variation in terms of their generating kernels.
\end{abstract}

\maketitle

\section{Introduction}

The definition of a variation of a real function was introduced by Camille Jordan as early as in 1890s, and since then numerous generalizations and extensions of this concept have been defined and thoroughly studied. One of them -- the Schramm variation introduced in~\cite{schramm} -- although technical in nature and seemingly artificial, is especially interesting for at least two reasons. First, in a single definition it embraces a plethora of other well-established notions of variation such as Young or Waterman, and allows to treat them in a unified manner. And, secondly, it has found applications in other fields of mathematics like, for example, approximation theory~\cite{MR4353508} and Henstock--Kurzweil integration theory~ \cite{MR3730314}.

Each year topics related to variation are gaining more and more popularity among mathematicians (see, for examples, the recent papers~\cites{wrobel,MR4162410,MR3542051} concerning the Schramm variation, or the general monographs~\cites{ABM,reinwand}).

One of the problems of the theory that stayed open for very long time was the full description of (relatively) compact subsets of spaces of functions of bounded variation.
Finally, in 2020 Bugajewski and Gulgowski, basing on a newly-introduced concept of a equi-variated set, established two compaction criteria the space $BV[0,1]$ of functions of bounded Jordan variation (see~\cite{BG20}*{Theorems~1 and~2}). This result was closely followed by analogous criteria in other spaces (see~\cite{G21}). In~\cite{GKM} we completed the study by establishing a compactness result in the space $\Phi BV[0,1]$ of functions of bounded Schramm variation.

Independently, Si and Xu in~\cite{SiXu} published a series of sufficient conditions for (relative) compactness in various spaces of functions of bounded variation. Their ideas directly followed the concept of the equi-variated set introduced in~\cite{BG20}. Unfortunately, although this approach fits the Jordan variation perfectly, it does not work for other variations, resulting in the loss of necessity of the considered conditions. It is worth noting here that in~\cite{GKM} we extended the notion of an equi-variated set. The method employed by Bugajewski and Gulgowski (and then repeated by Si and Xi) used a semi-norm related to a finite family of non-overlapping intervals. In our more general approach the considered family of intervals is still finite, but at the same time the intervals may overlap.

The aim of this paper is twofold. In the first part we revisit the compactness criterion in $\Phi BV[0,1]$. In particular, we show that some technical assumptions on the sequence $(\varphi_n)_{n \in \mathbb N}$ of
Young functions (like a $\Delta_2$-condition used in~\cite{SiXu}, or a monotonicity requirement used in~\cite{GKM}) are superfluous and can be removed.

In the second part of the paper we study linear integral operators acting between the spaces of functions of bounded Jordan and Schramm variations. Such operators are especially interesting, because they are the building blocks of non-linear Hammerstein or Volterra integral operators, which, in turn, play a key role in the theory of integral equations. In particular, we provide necessary and sufficient conditions for a linear integral operator to map the space $BV[0,1]$ into $\Phi BV[0,1]$ and to be continuous/compact, obtaining thus the full characterization of such operators in terms of their generating kernels. Those results extend the results obtained in~\cites{BCGS,PK, BGK} and complement the results of~\cite{AB} (see also the very recent monograph~\cite{reinwand}*{Section~4.3}).

\section{Preliminaries}

As usual $\mathbb N$ denotes the set of positive integers and $\mathbb R$ is the set of real numbers. By an \emph{interval} we understand a non-empty and convex subset of $\mathbb R$. In particular, we also admit \emph{degenerate} closed intervals consisting of a single element. If $I$ and $J$ are two intervals, we call them \emph{non-overlapping} when their intersection $I \cap J$ is either an empty set or a singleton. To simplify the notation, for a function $x$ and an interval $I:=[a,b]$, we write $x(I)$ for $x(b)-x(a)$. If $A\subseteq \mathbb R$, then by $\chi_A$ we denote the characteristic function of the set $A$.

In the paper we are going to use both the Lebesgue and the Riemann--Stieltjes integrals. Therefore, to avoid misunderstandings, we denote them by the symbols ``$(L)\int$'' and ``$(RS)\int$'', respectively. Since we will be working with the Schramm variation, we also need the notion of a Young function. A function $\varphi\colon [0,+\infty)\to[0,+\infty)$ is said to be a \emph{Young function} (or \emph{$\varphi$-function}) if it is convex and such that $\varphi(t)=0$ if and only if $t=0$. It can be checked without any difficulty that Young functions are continuous, strictly increasing and $\f(t)\to+\infty$ as $t\to+\infty$.

Finally, let us recall the two notions connected with compactness that will play a key role in our paper. A metric space $X$ is called \emph{precompact} if its completion is compact. Equivalently, $X$ is precompact if and only if every sequence in $X$ contains a Cauchy subsequence. It can be also proved that precompactness is equivalent to total boundedness. Let us recall that a metric space $X$ is called \emph{totally bounded} if for each $\varepsilon>0$ it has a finite \emph{$\varepsilon$-net}, that is, there exists a finite collection of points $a_1,\ldots,a_n \in X$ such that $X=\bigcup_{i=1}^n \ball_X(a_i,\varepsilon)$; here, $\ball_X(a_i,\varepsilon)$ denotes the closed ball in the metric space $X$ centered at $a_i$ and of radius $\varepsilon$. Let us also mention that a subset of a metric space $X$ is precompact if it is precompact as a metric space itself (with the metric inherited from $X$). On the other hand, a subset $A$ of a metric space $X$ is called \emph{relatively compact} in $X$ if its closure in $X$ is compact, or, equivalently, if every sequence of elements in $A$ contains a subsequence convergent in $X$. In general, relative compactness is a stronger notion than precompactness. However, if the underlying metric space $X$ is complete, those two notions coincide. For more information on precompact and relatively compact sets see, for example,~\cite{MV}*{Chapter 4}.

\section{Abstract precompactness criterion}
\label{sec:abstract_compactness_criterion}

Let $(E,\norm{\cdot}_E)$ be a normed space. Assume also that on $E$ we can define a family $\bigl\{\norm{\cdot}_i\bigr\}_{i \in I}$ of semi-norms which satisfies the following two conditions:
\begin{enumerate}[label=\textbf{\textup{(A\arabic*)}}]
 \item\label{i} $\norm{x}_E=\sup_{i \in I}\norm{x}_i$ for $x \in E$,
 \item\label{ii} for every $i,j \in I$ there is an index $k \in I$ such that $\norm{x}_i \leq \norm{x}_k$ and $\norm{x}_j \leq \norm{x}_k$ for $x \in E$.
\end{enumerate}
A non-empty subset $A$ of $E$ is called \emph{equinormed} (or more precisely, \emph{equinormed with respect to the family $\bigl\{\norm{\cdot}_i\bigr\}_{i \in I}$}) if for every $\varepsilon>0$ there exists an index $k \in I$ such that $\norm{x}_E\leq \varepsilon + \norm{x}_k$ for all $x \in A$.

In~\cite{GKM} the following abstract precompactness criterion was proven (see also~\cite{GKM2}*{Theorem~1} for its slight refinement).

\begin{theorem}\label{thm:compactness_ver2}
Let $(E,\norm{\cdot}_E)$ be a normed space equipped with a family $\bigl\{\norm{\cdot}_i\bigr\}_{i \in I}$ of semi-norms satisfying conditions~\ref{i} and~\ref{ii}. Furthermore, assume that 
\begin{enumerate}
 \item[\textup{\textbf{(A4)}}\footnotemark] for each bounded sequence $(x_n)_{n \in \mathbb N}$ of elements of $E$ there is a subsequence $(x_{n_k})_{k \in \mathbb N}$ which is Cauchy with respect to each semi-norm $\norm{\cdot}_i$.
\end{enumerate}
Then, a non-empty subset $A$ of $E$ is precompact if and only if it is bounded and the set $A-A$ is equinormed.
\end{theorem}
\footnotetext{We wanted the notation of this paper to be reader-friendly and consistent with the one we used in~\cite{GKM}. So, in particular, we decided to keep the original numbering of conditions and not to rename \textbf{(A4)} to \textbf{(A3)}.}

Looking at the definition of precompactness, it should not come as a surprise that we can replace the condition~\textbf{(A4)}, given in terms of sequences, with a condition formulated by means of $\varepsilon$-nets. In this way, we end up with a version of Theorem~\ref{thm:compactness_ver2}. For completeness, let us also add that in the statement of this result by $\ball_i(a,r)$ we understand the set $\dset{x \in E}{\norm{x-a}_i\leq r}$, where $\norm{\cdot}_i$ is a given semi-norm. 

\begin{theorem}\label{thm:compactness_ver3}
Let $(E,\norm{\cdot}_E)$ be a normed space equipped with a family $\bigl\{\norm{\cdot}_i\bigr\}_{i \in I}$ of semi-norms satisfying conditions~\ref{i} and~\ref{ii}. Furthermore, assume that
\begin{enumerate}[label=\textbf{\textup{(A\arabic*)}}]
\setcounter{enumi}{4}
 \item\label{iv1} for each bounded subset $B$ of $E$, any index $i\in I$ and $\varepsilon>0$ there exist finitely many points $x_1^i,\ldots,x_N^i \in B$ such that $B\subseteq \bigcup_{j=1}^N \ball_i(x_j^i,\varepsilon)$.
\end{enumerate}
Then, a non-empty subset $A$ of $E$ is precompact if and only if it is bounded and the set $A-A$ is equinormed.
\end{theorem}

\begin{proof}
The necessity part follows from the fact that in normed spaces precompactness carries over to algebraic sums and that precompact subsets are equinormed (see~\cite{GKM}*{Proposition~3.13}).

To prove the sufficiency, fix $\varepsilon>0$ and choose an index $k \in I$ such that $\norm{x-y}_E \leq \norm{x-y}_k + \frac{1}{2}\varepsilon$ for all $x,y\in A$. Moreover, since $A$ is bounded, by~\ref{iv1} there exists a finite collection of points $x_1^k,\ldots,x_N^k\in A$ such that $A \subseteq \bigcup_{j=1}^N \ball_i(x_j^k,\frac{1}{2}\varepsilon)$. To end the proof it is enough to note that this collection is a finite $\varepsilon$-net for $A$ in $E$.
\end{proof}

\section{Functions of bounded variation}

We are now going to recall the definition of the Schramm variation, with the approach that allows for degenerate intervals. 

Let $x$ be a real-valued function defined on $[0,1]$ and let $(\varphi_n)_{n \in \mathbb N}$ be a sequence of Young functions such that $\varphi_{k+1}(t)\leq \varphi_k(t)$ and $\sum_{n=1}^\infty \varphi_n(t)=+\infty$ for each $t >0$ and $k \in \mathbb N$. The (possibly infinite) quantity
\[
 \var_{\Phi} x :=\sup \sum_{n=1}^{k}\varphi_n\big(\abs{x(I_n)}\big),
\]
where the supremum is taken over all finite collections $I_1,\ldots,I_k$ of closed and non-overlapping subintervals of $[0,1]$, $k \in \mathbb N$, is called the \emph{Schramm variation} (or, the \emph{$\Phi$-variation}) of the function $x$ over $[0,1]$.

In the sequel we will also assume that the intervals $I_1,\ldots,I_k$ appearing in any collection are distinct. If all the intervals are non-degenerate, this requirement follows directly from the fact that they are non-overlapping. However, it may happen that a collection contains some degenerate intervals. And then, we simply remove the duplicates by changing the recurring single point to some other ,,unchosen'' real number. This procedure clearly does not change the sum appearing in the definition of $\var_{\Phi} x$. 

It is known that when we consider the sequence $(\varphi_n)_{n \in \mathbb N}$ of Young functions given by $\varphi_n(t)=t$ for $n \in \mathbb N$, then the Schramm variation reduces to the classical \emph{Jordan variation}. In this case instead of $\var_\Phi x$ we will simply write $\var x$.

A word of caution is needed here. In the original definition of the $\Phi$-variation provided by Schramm in his paper~\cite{schramm} it is not clear whether the considered intervals may be degenerate or not. Neither does the monograph on functions of bounded variation~\cite{ABM} shed any light on this problem. The thorough discussion of different approaches and the conclusion that all of them lead to the same result can be found in~\cite{GKM}*{Section~6}. 

A function $x \colon [0,1] \to \mathbb R$ is said to be of \emph{bounded Schramm variation} if there exists a number $\lambda>0$ such that $\var_\Phi(\lambda x)<+\infty$. For the Jordan variation this relation simplifies to $\var x < +\infty$. 

The linear space of all functions which are of bounded Schramm or Jordan variation will be denoted by $\Phi BV[0,1]$ and $BV[0,1]$, respectively. It is well-known that those spaces become Banach spaces when endowed with the norms $\norm{x}_\Phi: = \abs{x(0)}+\abs{x}_\Phi$, where 
\[
 \abs{x}_{\Phi}:=\inf\dset{\lambda>0}{\var_{\Phi}(\lambda^{-1}x)\leq 1},
\]
and $\norm{x}_{BV}:=\abs{x(0)}+\var x$, respectively (see, for example,~\cite{ABM}*{Proposition~2.44} and cf.~\cite{schramm}*{Theorem~2.3}). Observe also that the norm $\norm{\cdot}_{BV}$ is a special case of $\norm{\cdot}_{\Phi}$, as it it is not hard to prove that $\abs{x}_{\Phi}=\var x$, if we choose the Young functions $\varphi_n$ as $\varphi_n(t)=t$ for $n \in \mathbb N$ and $t \in [0,+\infty)$.

In the sequel we will also need a few properties of the functions of bounded Jordan variation (or, BV-functions for short). It is clear that such functions are bounded and $\norm{x}_{\infty}\leq \norm{x}_{BV}$ for $x \in BV[0,1]$. Furthermore, it can be shown that BV-functions can be represented as differences of two monotone functions, and hence they are Borel measurable. 

For a comprehensive treatment of functions of bounded variation in various senses we refer the reader to the monographs~\cite{ABM} and~\cite{reinwand}.

\section{Compactness criterion in $\Phi BV[a,b]$}

In this section we are going to revisit a compactness criterion in the space $\Phi BV[0,1]$ we established in~\cite{GKM}*{Section~6.2}. As we are going to apply the abstract result from Section~\ref{sec:abstract_compactness_criterion}, we need to define a suitable family of semi-norms on $\Phi BV[0,1]$.

Fix a sequence $(\varphi_n)_{n \in \mathbb N}$ of Young functions which for each $t>0$ and $k \in \mathbb N$ satisfies the following two conditions $\varphi_{k+1}(t)\leq \varphi_k(t)$ and $\sum_{n=1}^\infty \varphi_n(t)=+\infty$. Furthermore, let $\mathcal P$ be the family of all closed subintervals of $[0,1]$ and let $\mathcal J \subseteq \mathcal P$. For a function $x \colon [0,1] \to \mathbb R$ we set
\[
 V_{\mathcal J}(x):=\sup \sum_{n=1}^{k}\varphi_n\big(\abs{x(I_n)}\big),
\]
where the supremum is taken over all finite collections $I_1,\ldots,I_k$ of non-overlapping intervals such that $I_n \in \mathcal J$ for $n=1,\ldots,k$, $k \in \mathbb N$. If $x \in \Phi BV[0,1]$, this quantity is obviously finite. We also define
\begin{equation}\label{eq:semi-norm}
 \abs{x}_{\mathcal J}:=\inf\dset[\big]{\lambda>0}{V_{\mathcal J}(\tfrac{x}{\lambda})\leq 1}
\end{equation} 
and
\begin{equation}\label{eq:norm}
 \norm{x}_{\mathcal J}:=\abs{x(0)}+ \abs{x}_{\mathcal J}.
\end{equation} 
It turns out that for each family $\mathcal J$ formula~\eqref{eq:norm} defines a semi-norm on $\Phi BV[0,1]$ (cf.~\cite{ABM}*{Proposition~2.44}, \cite{GKM}*{Proposition~6.7}, or~\cite{schramm}*{Theorem~2.2}). Furthermore, it is clear that in the special case when $\mathcal J=\mathcal P$ we have $V_{\mathcal P}(x)=\var_{\Phi} x $, and consequently $\abs{x}_{\mathcal P}=\abs{x}_{\Phi}$ and $\norm{x}_{\mathcal P}=\norm{x}_{\Phi}$.

Now, we are in position to prove a refinement of~\cite{GKM}*{Proposition 6.10}. Since we used Helly's selection theorem for the $\Phi$-variation to prove the original result, we needed an additional assumption that the functions $t \mapsto \varphi_{n+1}(t)-\varphi_n(t)$ are non-increasing on $[0,+\infty)$ for $n \in \mathbb N$. (For a detailed discussion of this condition and the reason why it is important see~\cite{GKM}*{Remark~6.11}.) However, with a slight change in the approach this assumption becomes redundant. 

\begin{proposition}\label{prop:semi_norm_F}
Let $(\varphi_n)_{n \in \mathbb N}$ be a fixed sequence of Young functions which for each $t>0$ and $k \in \mathbb N$ satisfies the following two conditions $\varphi_{k+1}(t)\leq \varphi_k(t)$ and $\sum_{n=1}^\infty \varphi_n(t)=+\infty$. Then, the family of semi-norms $\{\norm{\cdot}_{\mathcal J}\}_{\mathcal J \subseteq \mathcal P}$ on $(\Phi BV[0,1], \norm{\cdot}_{\Phi})$ indexed by finite subfamilies $\mathcal J$ of $\mathcal P$, where $\norm{\cdot}_{\mathcal J}$ is given by the formula~\eqref{eq:norm}, satisfies the conditions~\ref{i}, \ref{ii} and~\ref{iv1}.      
\end{proposition}

\begin{proof}
The proof that the given family of semi-norms satisfies conditions~\ref{i} and \ref{ii} was given in \cite{GKM}*{Proposition~6.10}; the additional assumption on the monotonicity of the increments $\varphi_{n+1}-\varphi_n$ was  not used in that part. 

We need to prove that~\ref{iv1} is satisfied. Let us take any bounded subset $A$ of $\Phi BV[0,1]$, any finite collection $\mathcal J:=\{I_1,\ldots,I_N\}$ of closed subintervals of $[0,1]$ of the form $I_n:=[a_n,b_n]$ and an arbitrary $\varepsilon>0$. Define the linear projection $Q \colon \Phi BV[0,1] \to \mathbb R^{2N+1}$ corresponding to the family $\mathcal J$ by the formula 
\[
Q(x) = (x(0),x(b_1),x(a_1),\ldots,x(b_N),x(a_N));
\]
here, we endow the Euclidean space $\mathbb R^{2N+1}$ with the maximum norm $\norm{\cdot}_{\infty}$. By $\varphi_1^{-1} \colon [0,+\infty) \to [0,+\infty)$ let us denote the inverse function of $\varphi_1$; it exists, because Young functions are strictly increasing and map the interval $[0,+\infty)$ onto itself. As $\sup_{t \in [0,1]} \abs{x(t)} \leq (1+\varphi_1^{-1}(1))\norm{x}_{\Phi}$ for any $x \in \Phi BV[0,1]$ (see the proof of~\cite{ABM}*{Proposition~2.44~(d)} or~\cite{GKM}*{Remark~6.8}), the set $Q(A)$ is bounded (and, hence, relatively compact) in $\mathbb R^{2N+1}$. Therefore, it has a finite $\varepsilon (1+2\sigma)^{-1}$-net $Q(x_1),\ldots,Q(x_m)$, where $\sigma: = 1+ \sum_{n=1}^N \varphi_n(1)$. 

It remains to show that $x_1,\ldots,x_m \in A$ is an $\varepsilon$-net for $A$ with respect to the semi-norm $\norm{\cdot}_{\mathcal J}$. We will do it by proving that $\norm{x}_{\mathcal J} \leq (1+2\sigma)\norm{Q(x)}_{\infty}$ for $x \in \Phi BV[0,1]$. So, let us fix $x \in \Phi BV[0,1]$. Note that we may assume that $\norm{Q(x)}_\infty >0$. Otherwise, $V_{\mathcal J}(x)=0$. Consequently, $\norm{x}_{\mathcal J}=0$, and the estimate follows trivially. If $J_1,\ldots,J_k \in \mathcal J$ is any finite collection of non-overlapping intervals, then, by the convexity of the Young functions and the fact that $\varphi_n(0)=0$, we get
\[
 \sum_{n=1}^k \varphi_n\biggl( \frac{\abs{x(J_n)}}{2\sigma \norm{Q(x)}_\infty}\biggr) \leq \sum_{n=1}^k \frac{\abs{x(J_n)}}{2\sigma \norm{Q(x)}_\infty} \varphi_n(1) \leq 1.
\]     
Therefore, $\abs{x}_{\mathcal J} \leq 2\sigma \norm{Q(x)}_\infty$. Since we also have $\abs{x(0)}\leq \norm{Q(x)}_\infty$, the estimate $\norm{x}_{\mathcal J} \leq (1+2\sigma)\norm{Q(x)}_\infty$ follows. The proof is complete.   
\end{proof}

We end this section with a compactness criterion in $\Phi BV[0,1]$, which is a consequence of the previous result, Theorem~\ref{thm:compactness_ver3} and the fact that $\Phi BV[0,1]$ is a Banach space.

\begin{theorem}\label{thm:compactnes_PBV}
Let $(\varphi_n)_{n \in \mathbb N}$ be a fixed sequence of Young functions which for each $t>0$ and $k \in \mathbb N$ satisfies the following two conditions $\varphi_{k+1}(t)\leq \varphi_k(t)$ and $\sum_{n=1}^\infty \varphi_n(t)=+\infty$. For any bounded and non-empty subset $A$ of $(\Phi BV[0,1], \norm{\cdot}_\Phi)$ the following conditions are equivalent\textup:
\begin{enumerate}[label=\textup{(\roman*)}]
 \item $A$ is relatively compact in $(\Phi BV[0,1], \norm{\cdot}_\Phi)$,
 \item for every $\varepsilon>0$ there exists a finite family $\mathcal J$ of closed subintervals of $[0,1]$ such that $\abs{x-y}_\Phi \leq \varepsilon + \abs{x-y}_\mathcal J$ for all $x,y \in A$\textup;
\end{enumerate}
here the semi-norms $\abs{\cdot}_\Phi$ and $\abs{\cdot}_{\mathcal J}$ are defined by the formula~\eqref{eq:semi-norm}.
\end{theorem}

\begin{remark}
Theorem~\ref{thm:compactnes_PBV} refines~\cite{GKM}*{Theorem~6.12}. In comparison to that result, here we do not need to assume additionally that the functions $t \mapsto \varphi_{n+1}(t)-\varphi_n(t)$, where $n \in \mathbb N$, are non-increasing on $[0,+\infty)$. This also answers the question raised at the end of the paper~\cite{GKM}. 
\end{remark}

\begin{remark}
In a very recent paper~\cite{SiXu}, Yanan Si and Jingshi Xu proved a sufficient condition for a bounded set of $\Phi BV([0,1],E)$ to be relatively compact (Theorem~7.4); they used an extension of the notion of an equivariated set introduced by Bugajewski et al. in~\cite{BG20}. In comparison with our paper, Yanan Si and Jingshi Xu work with maps taking values in an arbitrary Banach space $E$. Moreover, they follow the approach presented in~\cite{ABM}*{Section~2.3} and do not assume in the definition of the Schramm variation that the sequence $(\varphi_n)_{n \in \mathbb N}$ of Young functions is pointwise non-increasing, that is, $\varphi_{n+1}(t)\leq \varphi_n(t)$ for all $n \in \mathbb N$ and $t \in [0,+\infty)$. On the other hand, they require the sequence $(\varphi_n)_{n \in \mathbb N}$ to satisfy a $\Delta_2$-type condition. Namely, it is assumed that there exists a constant $M> 0$ such that $\varphi_n(2t) \leq M\varphi_n(t)$ for all $t \in [0,+\infty)$ and $n \in \mathbb N$. Although very natural in the theory of modular spaces, $\Delta_2$-conditions can be also restrictive.    
\end{remark}

\section{Integral operators}

In this section we will continue our investigation of functions of Schramm bounded variation. This time, however, we will be interested in linear integral operators. More precisely, we are going to provide necessary and sufficient conditions on the kernel under which the integral operator, generated by this kernel, maps the space $BV[0,1]$ into $\Phi BV[0,1]$ and is continuous and/or compact. We will follow the approach introduced in~\cite{BGK} and developed in~\cite{PK}. The proofs, however, will be much more technical and demanding. That is also the reason we decided to provide them in full detail.

As our approach is based on rewriting the Lebesgue integral as the Riemann--Stieltjes integral, let us recall now three facts about the Riemann--Stieltjes integration. The first one, in much simpler form, is usually taught during the Analysis 101 course. 

\begin{proposition}[see~\cite{PK}*{Proposition~2} and cf.~\cite{Lojasiewicz}*{Theorem~7.4.10}]\label{prop:reduction}
Let $f \colon [0,1] \to \mathbb R$ be a Lebesgue integrable function and let $g \colon [0,1] \to \mathbb R$ be of bounded Jordan variation. Then,
\begin{equation*}\label{eq:L_RS}
 (L) \int_0^1 f(t)g(t) \textup dt = (RS) \int_0^1 g(t) \textup dF(t),
\end{equation*}
where $F$ is the primitive of $f$, that is, $F(t)=(L)\int_0^t f(s)\textup ds$ for $t \in [0,1]$.
\end{proposition}

Another two results which will come in handy are the integration-by-parts formula and Jensen's inequality for the Riemann--Stieltjes integral. We do not state those theorems in their full generality, but in a form which we will need in a sequel. Moreover, we decided to provide the proof of Proposition~\ref{prop:jensen} for readers' convenience as we could not find it in any book.

\begin{proposition}[see~\cite{Lojasiewicz}*{Theorem~1.6.7}]\label{prop:parts}
If $f,g \colon [0,1] \to \mathbb R$ are functions of bounded Jordan variation which do not have common points of discontinuity, then
\[
 (RS) \int_0^1 f(t) \textup dg(t) + (RS) \int_0^1 g(t) \textup df(t) = f(1)g(1) - f(0)g(0).
\]
\end{proposition}

\begin{proposition}\label{prop:jensen}
Let $\varphi \colon [0,+\infty) \to [0,+\infty)$ be a Young function. If $f \colon [0,1] \to \mathbb [0,+\infty)$ is a continuous function and $g \colon [0,1] \to [0,1]$ is a non-decreasing function, then
\begin{equation}\label{eq:jensen}
 \varphi\biggl((RS)\int_0^1 f(t) \textup dg(t)\biggr) \leq (RS)\int_0^1 \varphi(f(t))\textup dg(t).
\end{equation}
\end{proposition}

\begin{proof}
Fix an arbitrary positive integer $n$. In view of the convexity of the function $\varphi$ and the facts that $\varphi(0)=0$ and $g(t) \in [0,1]$ for $t \in [0,1]$, we have
\[
 \varphi\biggl(\sum_{j=1}^n f(\tfrac{j}{n})[g(\tfrac{j}{n})-g(\tfrac{j-1}{n})]\biggr) \leq \sum_{j=1}^n \varphi(f(\tfrac{j}{n}))[g(\tfrac{j}{n})-g(\tfrac{j-1}{n})].
\] 
Since both the integrals in~\eqref{eq:jensen} exist and are finite (cf.~\cite{Lojasiewicz}*{Theorem~1.5.5}), to obtain the estimate it suffices now to pass to the limit with $n \to +\infty$ in the above inequality.
\end{proof}

Now, let us move to the main part of this section. Given a kernel $k \colon [0,1] \times [0,1] \to \mathbb R$ let us consider the linear integral operator of the form
\begin{equation}\label{eq:opK}
 (Kx)(t) = (L)\int_0^1 k(t,s) x(s) \textup ds, \quad t \in [0,1],
\end{equation}
where $x \in BV[0,1]$. Clearly, the formula~\eqref{eq:opK} makes sense for any $x \in BV[0,1]$ only if the function $k$ satisfies the following condition
\begin{enumerate}[leftmargin=35pt, itemsep=2pt, label=\textup{\textbf{(H\arabic*)}}] 
 \item\label{H1} for every $t \in [0,1]$ the function $s \mapsto k(t,s)$ is Lebesgue integrable on $[0,1]$.
\end{enumerate}

It may come as a surprise that only with the definition of the integral operator and a few auxiliary results from mathematical analysis at hand we are ready to prove necessary and sufficient conditions for the operator $K$ to map $BV[0,1]$ into $\Phi BV[0,1]$ continuously.   

\begin{theorem}\label{thm:K_bV_PBV}
Let $(\varphi_n)_{n \in \mathbb N}$ be a fixed sequence of Young functions such that $\varphi_{j+1}(t)\leq \varphi_j(t)$ and $\sum_{n=1}^\infty \varphi_n(t)=+\infty$ for each $t >0$ and $j \in \mathbb N$. Moreover, let $K$ be the integral operator given by~\eqref{eq:opK} whose kernel $k \colon [0,1] \times [0,1] \to \mathbb R$ satisfies condition~\ref{H1}. The operator $K$ maps the space $BV[0,1]$ into $\Phi BV[0,1]$ and is continuous if and only if 
\begin{enumerate}[leftmargin=35pt, itemsep=2pt, label=\textup{\textbf{(H\arabic*)}}, ]
\setcounter{enumi}{1} 
 \item\label{H2} there exists a constant $\mu>0$ such that for every $\xi \in [0,1]$ we have
\[
\var_{\Phi}\biggl(\mu \cdot (L)\int_0^{\xi} k(\cdot,s)\textup ds\biggr) \leq 1.
\]
\end{enumerate}
\end{theorem}

\begin{proof}
First, let us assume that the operator $K$ maps the space $BV[0,1]$ into $\Phi BV[0,1]$ and is continuous. Then, there exits a constant $M>0$ such that $\norm{K(\chi_{[0,\xi]})}_{\Phi} \leq M\norm{\chi_{[0,\xi]}}_{BV}\leq 2M$ for each $\xi \in [0,1]$. Therefore, for every $\xi \in [0,1]$ we have
\[
 \var_{\Phi}\biggl(\mu \cdot (L)\int_0^{\xi}  k(\cdot,s)\textup ds\biggr) \leq 1.
\]
with, for example, $\mu:=(2M)^{-1}$ (cf.~\cite{schramm}*{Lemma~2.1~(i)}). This ends the proof of the necessity.

To show the sufficiency, let us assume that $k$ satisfies~\ref{H2}. If $x \in BV[0,1]$ is a constant function, then
\[
 (Kx)(t)=x(0)\cdot (L)\int_0^1 k(t,s)\textup ds
\]
for $t \in [0,1]$, and hence $Kx\in \Phi BV[0,1]$. Furthermore, 
\begin{align}\label{eq:continuity_constant_x}
\begin{split}
 \norm{Kx}_{\Phi} & \leq \abs{x(0)}\cdot \abs[\bigg]{(L)\int_0^1 k(0,s)\textup ds}+\abs{x(0)}\cdot \abs[\bigg]{(L)\int_0^1 k(t,s)\textup ds}_\Phi\\
 &\leq \biggl((L)\int_0^1 \abs{k(0,s)}\textup ds + \mu^{-1}\biggr) \cdot \norm{x}_{BV}.
\end{split}
\end{align}

Now, let us assume that $x \in BV[0,1]$ is non-constant; this ensures that $\var x>0$. Note that for $t \in [0,1]$ in view of Propositions~\ref{prop:reduction} and~\ref{prop:parts} we have
\begin{equation}\label{eq:representation_of_Kx}
 (Kx)(t)= x(1)\cdot (L)\int_0^1 k(t,s)\textup ds - (RS)\int_0^1 \biggl((L)\int_0^\xi k(t,s)\textup ds\biggr) \textup dx(\xi).
\end{equation}
To simplify the calculations, we will extend the notation we used in the previous sections and given an interval $J:=[a,b]$ we will write $k(J,s)$ for $k(b,s)-k(a,s)$. If $\lambda>0$ and $I_1,\ldots,I_n$ is a finite collection of closed and non-overlapping subintervals of $[0,1]$, then
\begin{align*}
 \sum_{j=1}^n \varphi_j\biggl(\frac{\abs{(Kx)(I_j)}}{\lambda}\biggr) \leq \sum_{j=1}^n \varphi_j\biggl(\frac{\abs{x(1)}}{\lambda}\cdot \abs[\bigg]{(L)\int_0^1 k(I_j,s)\textup ds} + \frac{1}{\lambda}\cdot (RS)\int_0^1 \abs[\bigg]{(L)\int_0^\xi k(I_j,s)\textup ds} \textup dv_x(\xi)\biggr), 
\end{align*}
where $v_x$ denotes the variation function given by $v_x(0)=0$ and $v_x(t)=\var(x,[0,t])$ for $t \in (0,1]$, and the symbol $\var(x,[0,t])$ stands for the Jordan variation of the function $x$ over the interval $[0,t]$ -- the definition is analogous to the one for the interval $[0,1]$. In the above inequality we also used the classical estimate for the Riemann--Stieltjes integral which states that for a continuous function $f\colon [0,1] \to \mathbb R$ and a function of bounded Jordan variation $g \colon [0,1] \to \mathbb R$ we have
\[
 \abs[\bigg]{(RS)\int_0^1 f(t)dg(t)}\leq (RS)\int_0^1 \abs{f(t)}dv_g(t)
\] 
(see~\cite{ABM}*{Theorem~4.20}). Setting $\lambda:=2\mu^{-1}\norm{x}_{BV}$ and $\overline{v}_x(\xi):=v_x(\xi)/\var x$, as the function $v_x$ is non-decreasing, by Proposition~\ref{prop:jensen} and properties of $\varphi_j$'s, we obtain
\begin{align*}
& \sum_{j=1}^n \varphi_j\biggl(\frac{\abs{(Kx)(I_j)}}{\lambda}\biggr)\\
&\quad \leq \frac{1}{2}\sum_{j=1}^n \varphi_j \biggl(\abs[\bigg]{\mu \cdot (L)\int_0^1 k(I_j,s)\textup ds}\biggr) + \frac{1}{2}\sum_{j=1}^n\varphi_j \biggl((RS)\int_0^1 \abs[\bigg]{\mu \cdot (L)\int_0^\xi k(I_j,s)\textup ds} \textup d\overline{v}_x(\xi)\biggr)\\
& \quad \leq \frac{1}{2}\var_\Phi\biggl(\mu \cdot (L)\int_0^1 k(\cdot,s)\textup ds\biggr) + \frac{1}{2} (RS)\int_0^1 \sum_{j=1}^n \varphi_j \biggl(\abs[\bigg]{\mu\cdot (L)\int_0^\xi k(I_j,s)\textup ds}\biggr) \textup d\overline{v}_x(\xi).
\end{align*}
As
\[
\sum_{j=1}^n \varphi_j \biggl(\abs[\bigg]{\mu\cdot (L)\int_0^\xi k(I_j,s)\textup ds}\biggr) \leq \var_\Phi\biggl(\mu \cdot (L)\int_0^\xi k(\cdot,s)\textup ds\biggr) \leq 1
\]
for $\xi \in [0,1]$, we get
\[
  \sum_{j=1}^n \varphi_j\biggl(\frac{\abs{(Kx)(I_j)}}{\lambda}\biggr)
  \leq \frac{1}{2}+ \frac{1}{2}(RS)\int_0^1 1 \textup d\overline{v}_x(\xi) = 1.
\]
This shows that $Kx \in \Phi BV[0,1]$ and that $\abs{Kx}_{\Phi} \leq 2\mu^{-1}\norm{x}_{BV}$. Hence, 
\[
\norm{Kx}_{\Phi} \leq \norm{x}_{\infty} \cdot (L)\int_0^1 \abs{k(0,s)}\textup ds + 2\mu^{-1}\norm{x}_{BV}.
\]
Together with~\eqref{eq:continuity_constant_x} and the fact that $\norm{x}_\infty \leq \norm{x}_{BV}$ this proves that $K$ is continuous and $\norm{Kx}_{\Phi}\leq M\norm{x}_{BV}$ with
\[
 M:=(L)\int_0^1 \abs{k(0,s)}\textup ds + 2\mu^{-1}.
\]
The proof is complete.
\end{proof}

\begin{remark}
Theorem~\ref{thm:K_bV_PBV} extends analogous results established in~\cites{BGK,BCGS} for linear integral operators acting into the spaces $BV[0,1]$, $BV_p[0,1]$ and $\Lambda BV[0,1]$ of functions of bounded Jordan, Wiener and Waterman variation, respectively (cf. also~\cite{reinwand}*{Section~4.3}).  
\end{remark}

In 2016 Bugajewski et al., studying linear integral operators acting in $BV[0,1]$, provided an example of an operator with a kernel satisfying condition~\ref{H1} and a condition analogous to~\ref{H2} for functions of Jordan bounded variation which mapped the space $BV[0,1]$ into itself and was continuous but not compact (see~\cite{BGK}*{Example~3}). Consequently, in general, condition~\ref{H2} cannot guarantee compactness of a linear integral operator $K \colon BV[0,1] \to \Phi BV[0,1]$. To find appropriate necessary and sufficient conditions for compactness of such a map, we will follow the approach introduced in~\cite{PK}. Although, the proof of the following technical result is almost identical to the proof of~\cite{PK}*{Proposition~4}, we provide it for completeness. This proposition may be of independent interest, as it shows that when studying compactness of integral operators acting into various spaces of functions of bounded variation, we need to focus our attention only on those bounded sequences that are pointwise convergent to $0$. 

\begin{proposition}\label{prop:compactness_equiv}
Let $(\varphi_n)_{n \in \mathbb N}$ be a fixed sequence of Young functions such that $\varphi_{j+1}(t)\leq \varphi_j(t)$ and $\sum_{n=1}^\infty \varphi_n(t)=+\infty$ for each $t >0$ and $j \in \mathbb N$. Moreover, let $k\colon [0,1] \times [0,1] \to \mathbb R$ satisfy conditions~\ref{H1} and~\ref{H2} and let $K \colon BV[0,1] \to \Phi BV[0,1]$ be the linear integral operator given by~\eqref{eq:opK}. Then, the following conditions are equivalent\textup:
\begin{enumerate}[label=\textup{(\roman*)}]
 \item\label{prop:a} the operator $K$ is compact, that is, the image $K(\ball_{BV}(0,1))$ of the closed unit ball in $BV[0,1]$ is a relatively compact subset of $\Phi BV[0,1]$\textup,
 
  \item\label{prop:b} for every sequence $(x_n)_{n \in \mathbb N}$ of elements of $\ball_{BV}(0,1)$ which is pointwise convergent on $[0,1]$ to a function $x \colon [0,1] \to \mathbb R$, we have $\lim_{n \to \infty}\norm{Kx_n - Kx}_{\Phi}=0$\textup,
  
  \item\label{prop:c} for every sequence $(x_n)_{n \in \mathbb N}$ of elements of $\ball_{BV}(0,1)$ which is pointwise convergent on $[0,1]$ to the zero function, we have $\lim_{n \to \infty}\norm{Kx_n}_{\Phi}=0$.
\end{enumerate}
\end{proposition}

\begin{remark}
Of course, in the statement of Proposition~\ref{prop:compactness_equiv} we may replace the closed unit ball in $BV[0,1]$ with any closed ball centered at zero.
\end{remark}

\begin{proof}[Proof of Proposition~\ref{prop:compactness_equiv}]
$\ref{prop:a}\Rightarrow\ref{prop:b}$ Let $(x_n)_{n\in\mathbb N}$ be an arbitrary sequence with terms lying in the unit ball of $BV[0,1]$ which is pointwise convergent on $[0,1]$ to a function $x\colon [0,1] \to \mathbb R$. Then, $x$ is also of bounded Jordan variation and $x \in \ball_{BV}(0,1)$ -- see~\cite{Lojasiewicz}*{Theorem~1.3.5}. This, together with Theorem~\ref{thm:K_bV_PBV}, implies that the function $Kx$ is well-defined and belongs to $\Phi BV[0,1]$. Now, it remains to show that $\norm{Kx_n - Kx}_{\Phi} \to 0$ as $n \to +\infty$. First, by the Lebesgue dominated convergence theorem and assumption~\ref{H1}, we deduce that $(Kx_n)(t) \to (Kx)(t)$ pointwise on $[0,1]$. Next, we consider an arbitrary subsequence $(Kx_{n_m})_{m \in \mathbb N}$ of $(Kx_{n})_{n \in \mathbb N}$. Since the operator $K$ is compact, $(Kx_{n_m})_{m \in \mathbb N}$ contains yet another subsequence $(Kx_{n_{m_l}})_{l \in \mathbb N}$ which is convergent in the $\norm{\cdot}_\Phi$-norm to a function $y$ of Schramm bounded variation on $[0,1]$. Recall that for every $z \in \Phi BV[0,1]$ we have $\norm{z}_\infty \leq (1+\varphi_1^{-1}(1))\norm{z}_{\Phi}$, where $\varphi_1^{-1}$ denotes the inverse function of $\varphi_1$. Hence, the norm convergence in $\Phi BV[0,1]$ implies the pointwise convergence. So, $y(t)=\lim_{l \to +\infty} (Kx_{n_{m_l}})(t)$ for every $t \in [0,1]$. This means that $y=Kx$ and proves that $Kx_n \to Kx$ in $\Phi BV[0,1]$.

$\ref{prop:b}\Rightarrow\ref{prop:c}$ The implication is obvious.

$\ref{prop:c}\Rightarrow\ref{prop:a}$ Let $(x_n)_{n\in \mathbb N}$ be an arbitrary sequence of elements of the closed unit ball in $BV[0,1]$. Using the classical Helly's selection theorem (see, for example,~\cite{C}*{Theorem~13.16}), we deduce that there is a subsequence $(x_{n_m})_{m \in \mathbb N}$ of $(x_n)_{n \in \mathbb N}$ which is pointwise convergent on $[0,1]$ to some $x \in \ball_{BV}(0,1)$. Setting $y_m:=\frac{1}{2}(x_{n_m}-x)$, we obtain a sequence of elements of the closed unit ball in $BV[0,1]$ which is pointwise convergent on $[0,1]$ to $0$. Therefore, by our assumption, we get $\norm{Ky_m}_{\Phi}\to 0$ as $m \to +\infty$, but this implies that $\norm{Kx_{n_m} - Kx}_{\Phi} \to 0$ with $m \to +\infty$ and completes the proof.
\end{proof}

Equipped with Proposition~\ref{prop:compactness_equiv} we are ready to prove our second main theorem of this section.

\begin{theorem}\label{thm:K_BV_PBV_compact}
Let $(\varphi_n)_{n \in \mathbb N}$ be a fixed sequence of Young functions such that $\varphi_{j+1}(t)\leq \varphi_j(t)$ and $\sum_{n=1}^\infty \varphi_n(t)=+\infty$ for each $t >0$ and $j \in \mathbb N$. Moreover, let $k \colon [0,1] \times [0,1] \to \mathbb R$ be a kernel satisfying condition~\ref{H1} and let $K$ be the linear integral operator given by~\eqref{eq:opK}. The operator $K$ maps the space $BV[0,1]$ into $\Phi BV[0,1]$ and is compact if and only if
\begin{enumerate}[leftmargin=35pt, itemsep=2pt, label=\textup{\textbf{(H\arabic*)}}] 
\setcounter{enumi}{2}
\item\label{H3} for every $\varepsilon>0$ there exists $\delta:=\delta(\varepsilon)>0$ such that
\[
  \var_{\Phi}\biggl(\varepsilon^{-1} \cdot (L)\int_a^b  k(\cdot,s)\textup ds\biggr) \leq 1
\]
for any subinterval $[a,b]$ of $[0,1]$ of length not exceeding $\delta$.
\end{enumerate}
\end{theorem}

\begin{remark}\label{rem:h3_implies_h2}
Observe that if the kernel $k \colon [0,1] \times [0,1] \to \mathbb R$ satisfies condition~\ref{H1}, then~\ref{H3} implies \ref{H2}. To show our claim let $\delta:=\delta(1)>0$ be as in~\ref{H3}. Moreover, let us take a positive integer $n$ such that $\frac{1}{n} \leq \delta$ and set $\mu:=\frac{1}{n}$. Now, we fix $\xi \in [0,1]$ and we divide the interval $(0,1]$ into $n$ subintervals $(\frac{i-1}{n},\frac{i}{n}]$ of length $\frac{1}{n}$; here, $i=1,\ldots,n$. We also choose a number $m\in\{1,\ldots,n\}$ so that $\xi \in (\frac{m-1}{n},\frac{m}{n}]$. Then, using the convexity of $\var_\Phi$ and the fact that $\var_\Phi 0=0$ we have
\begin{align*}
 &\var_\Phi\biggl(\frac{1}{n}\cdot (L)\int_0^\xi k(\cdot,s)\textup ds\biggr)\\
 &\qquad \leq \sum_{i=1}^{m-1} \frac{1}{n} \var_\Phi\biggl((L)\int_{\frac{i-1}{n}}^{\frac{i}{n}} k(\cdot,s)\textup ds\biggr) + \frac{1}{n}\var_\Phi \biggl((L)\int_{\frac{m-1}{n}}^{\xi} k(\cdot,s)\textup ds\biggr) \leq \frac{m}{n}\leq 1.
\end{align*}
(If $m$ happens to be equal to $1$, by definition, we set the sum in the above estimate to be zero.) Thus, condition~\ref{H2} is satisfied with $\mu=\frac{1}{n}$.

Conditions~\ref{H2} and~\ref{H3} cannot be equivalent. As already mentioned, an appropriate example in the case of the Jordan variation was given by Bugajewski et al. in~\cite{BGK}*{Example~3}. (A word of warning: the condition \textbf{(H3)} in that paper is different from the condition \textbf{(H3)} we use here.)
\end{remark}

\begin{proof}[Proof of Theorem~\ref{thm:K_BV_PBV_compact}]
Let us assume that the integral operator $K$ given by~\eqref{eq:opK} maps the space $BV[0,1]$ into $\Phi BV[0,1]$ and is compact. Then, the image $K(\ball_{BV}(0,2))$ of the closed ball in $BV[0,1]$ is a relatively compact subset of $\Phi BV[0,1]$. In view of Theorem~\ref{thm:compactnes_PBV}  this, in turn, implies that given $\varepsilon>0$ there is a family $\mathcal J:=\{I_1,\ldots,I_n\}$ of closed (and \emph{not} necessarily non-overlapping) subintervals of $[0,1]$ of the form $I_j:=[a_j,b_j]$ such that $\abs{K\chi_{(a,b)}}_{\Phi} \leq \frac{1}{2}\varepsilon + \abs{K\chi_{(a,b)}}_{\mathcal J}$ for every interval $(a,b) \subseteq [0,1]$; here, the semi-norms $\abs{\cdot}_\Phi$ and $\abs{\cdot}_{\mathcal J}$ are defined by the formula~\eqref{eq:semi-norm}. Since by~\ref{H1} the function $s \mapsto k(t,s)$ is Lebesgue integrable for every $t\in [0,1]$, there exists a positive number $\delta$ such that for any $j \in \{1,\ldots,n\}$, any $t_j \in \{a_j,b_j\}$ and any Lebesgue measurable set $A \subseteq [0,1]$ of measure not exceeding $\delta$ we have
\[
 (L)\int_A \abs{k(t_j,s)}\textup ds \leq \tfrac{1}{8}\varepsilon \varphi_1^{-1}(\tfrac{1}{n});
\]  
by $\varphi_1^{-1}$ we denote the inverse function of $\varphi_1$. Now, let $(a,b) \subseteq [0,1]$ be an arbitrary (but fixed) interval such that $b-a\leq \delta$ and let $J_1,\ldots,J_m$ be a finite collection of non-overlapping intervals so that $J_l:=[c_l,d_l] \in \mathcal J$ for $l=1,\ldots,m$. Then,
\begin{align*}
&\sum_{l=1}^m \varphi_l\bigl(4\varepsilon^{-1}\abs[\big]{(K\chi_{(a,b)})(d_l)-(K\chi_{(a,b)})(c_l)}\bigr)\\
&\qquad \leq \sum_{l=1}^m \varphi_l\biggl (4\varepsilon^{-1} (L)\int_a^b \abs{k(d_l,s)-k(c_l,s)}\textup ds\biggr)\\
&\qquad \leq \frac{1}{2} \sum_{l=1}^m \varphi_l\biggl (8\varepsilon^{-1} (L)\int_a^b \abs{k(d_l,s)}\textup ds\biggr) + \frac{1}{2}\sum_{l=1}^m \varphi_l\biggl (8\varepsilon^{-1} (L)\int_a^b \abs{k(c_l,s)}\textup ds\biggr)\\
&\qquad \leq \frac{1}{2} \sum_{l=1}^m \varphi_l(\varphi_1^{-1}(\tfrac{1}{n})) + \frac{1}{2} \sum_{l=1}^m \varphi_l(\varphi_1^{-1}(\tfrac{1}{n}))\\
&\qquad \leq \sum_{l=1}^m \varphi_1(\varphi_1^{-1}(\tfrac{1}{n}))= \tfrac{m}{n}\leq 1,
\end{align*}
because $\dset{c_l,d_l}{l=1,\ldots,m}\subseteq \dset{a_j,b_j}{j=1,\ldots,n}$ and the cardinality of $\mathcal J$ is not smaller than $m$. Thus, $V_{\mathcal J}(4\varepsilon^{-1}K\chi_{(a,b)})\leq 1$, which implies that $\abs{K\chi_{(a,b)}}_{\Phi} \leq \tfrac{1}{2}\varepsilon + \abs{K\chi_{(a,b)}}_{\mathcal J} \leq  \frac{1}{2}\varepsilon + \frac{1}{4}\varepsilon<\varepsilon$. And so, 
\[
\var_{\Phi}\biggl(\varepsilon^{-1} \cdot (L)\int_a^b  k(t,s)\textup ds\biggr) = \var_{\Phi}(\varepsilon^{-1} K\chi_{(a,b)}) \leq 1.
\]
This means that the kernel $k$ satisfies condition~\ref{H3}.

Now, we prove the opposite implication. This time we assume that $k$ satisfies condition~\ref{H3} (and, hence, also~\ref{H2}). We are going to show that the integral operator $K$ is compact. To this end we will use Proposition~\ref{prop:compactness_equiv}. So, let us consider a sequence $(x_v)_{v \in \mathbb N}$ of elements of the closed unit ball $\ball_{BV}(0,1)$ which is pointwise convergent on $[0,1]$ to the zero function. Let us also fix $\varepsilon>0$ and consider a fixed partition  $\Xi\!: 0=\xi_0<\xi_1<\cdots<\xi_n=1$ of the interval $[0,1]$ such that $\max_{1\leq i \leq n}\abs{\xi_i - \xi_{i-1}}\leq \delta$, where $\delta:=\delta(\frac{1}{4}\varepsilon)>0$ is chosen as in~\ref{H3}. Furthermore, choose $v_0 \in \mathbb N$ so that for every $v\geq v_0$ we have 
\[
 \abs{x_v(1)}\leq \tfrac{1}{4}\varepsilon\mu \quad \text{and}\quad \sum_{i=1}^n \abs{x_v(\xi_i) - x_v(\xi_{i-1})} \leq \tfrac{1}{4}\varepsilon\mu,
\]
where $\mu$ is the constant appearing in~\ref{H2} (cf.~Remark~\ref{rem:h3_implies_h2}).

We now rewrite the function $Kx_v$ in such a way that will allow us to estimate its Schramm variation. So, let us fix $v \geq v_0$ and let $I_1,\ldots,I_m$ be a finite collection of closed and non-overlapping subintervals of $[0,1]$. As in the proof of Theorem~\ref{thm:K_bV_PBV}, given an interval $J:=[a,b]$ we will write $k(J,s)$ for $k(b,s)-k(a,s)$. Using assumption~\ref{H1} we infer that for each $j\in \{1,\ldots,m\}$ the function $\xi \mapsto (L)\int_0^{\xi} k(I_j,s)\textup ds$, where $\xi \in [0,1]$, is absolutely continuous, and so the Riemann--Stieltjes integral
\[
(RS) \int_0^1 \biggl((L)\int_0^{\xi} k(I_j,s)\textup ds \biggr) \textup dx_v(\xi)
\]
exists (cf.~\cite{Lojasiewicz}*{Theorem~1.5.5}). This, in turn, means that there is a partition 
 $\overline{\Xi}\!: 0=\overline{\xi}_0<\overline{\xi}_1<\cdots<\overline{\xi}_p=1$, which is a refinement of $\Xi$, such that
\begin{equation}\label{eq:estimate_of_RS_sum}
\sum_{j=1}^m\abs[\bigg]{\overline{S}_j^v - (RS) \int_0^1 \biggl((L)\int_0^{\xi}k(I_j,s)\textup ds \biggr)\textup dx_v(\xi)}\leq \frac{\varepsilon}{4(1+\varphi_1(1))}, 
\end{equation}
where $\overline{S}_j^v$ is the approximating Riemann--Stieltjes sum, that is, 
\[
  \overline{S}_j^v :=\sum_{l=1}^p\biggl( (L)\int_0^{\overline{\xi}_l} k(I_j,s)\textup ds\biggr)[x_v(\overline{\xi}_l)-x_v(\overline{\xi}_{l-1})].
\]
(As in the simple version of this proof -- see~\cite{PK}*{Theorem~2} -- the partition $\overline{\Xi}$ depends on the choice of $v$. However, since we assume that $v$ is fixed at the moment, we decided not to complicate the notation and not to write $\overline{\Xi}^v$ or $\overline{\xi}_l^v$.) Setting
\[
  S_j^v:=\sum_{i=1}^n \biggl( (L)\int_0^{\xi_i} k(I_j,s)\textup ds\biggr)[x_v(\xi_i)-x_v(\xi_{i-1})]
\]
and noting that $\xi_i=\overline{\xi}_{\alpha_i}$, where $0=\alpha_0<\cdots<\alpha_n=p$, for $i \in\{0,\ldots,n\}$, we have
\[
 S_j^v=\sum_{i=1}^n \sum_{l=\alpha_{i-1}+1}^{\alpha_i}\biggl( (L)\int_0^{\xi_i} k(I_j,s)\textup ds\biggr)[x_v(\overline{\xi}_l)-x_v(\overline{\xi}_{l-1})] 
\]
and
\[
 \overline{S}_j^v =\sum_{i=1}^n \sum_{l=\alpha_{i-1}+1}^{\alpha_i}\biggl( (L)\int_0^{\overline{\xi}_l} k(I_j,s)\textup ds\biggr)[x_v(\overline{\xi}_l)-x_v(\overline{\xi}_{l-1})].
\]
Finally, rewriting $Kx_v$ as in~\eqref{eq:representation_of_Kx} and using the properties of the Young functions $\varphi_j$ we have
\begin{align*}
& \sum_{j=1}^m \varphi_j\bigl(\varepsilon^{-1}\abs{(Kx_v)(I_j)}\bigr)\\
&\quad \leq \frac{1}{4}\sum_{j=1}^m \varphi_j\biggl(4\varepsilon^{-1} \abs{x_v(1)}\cdot \abs[\bigg]{(L)\int_0^1 k(I_j,s)\textup ds}\biggr)\\
&\qquad{}+ \frac{1}{4}\sum_{j=1}^m \varphi_j\biggl(4\varepsilon^{-1} \abs[\bigg]{(RS) \int_0^1 \biggl((L)\int_0^{\xi}k(I_j,s)\textup ds \biggr)\textup dx_v(\xi)-\overline{S}_j^v}\biggr)\\
&\qquad{}+  \frac{1}{4}\sum_{j=1}^m \varphi_j(4\varepsilon^{-1}\abs{\overline{S}_j^v - S_j^v}) +  \frac{1}{4}\sum_{j=1}^m \varphi_j(4\varepsilon^{-1}\abs{S_j^v}).
\end{align*}
Now, we estimate each of those four summands separately. We begin with the upper one. We have
\begin{align*}
& \sum_{j=1}^m \varphi_j\biggl(4\varepsilon^{-1} \abs{x_v(1)}\cdot \abs[\bigg]{(L)\int_0^1 k(I_j,s)\textup ds}\biggr)\\
&\qquad \leq \sum_{j=1}^m \varphi_j\biggl(\mu\abs[\bigg]{(L)\int_0^1 k(I_j,s)\textup ds}\biggr) \leq \var_\Phi \biggl(\mu \cdot (L)\int_0^1 k(\cdot,s)\textup ds\biggr)\leq 1.
\end{align*}
As regards the next one, since $1/(1+\varphi_1(1))\leq 1$, we obtain
\begin{align*}
&\sum_{j=1}^m \varphi_j\biggl(4\varepsilon^{-1} \abs[\bigg]{(RS) \int_0^1 \biggl((L)\int_0^{\xi} k(I_j,s)\textup ds \biggr)\textup dx_v(\xi)-\overline{S}_j^v}\biggr)\\
&\qquad \leq \sum_{j=1}^m \varphi_j(1) \cdot 4\varepsilon^{-1} \abs[\bigg]{(RS) \int_0^1 \biggl((L)\int_0^{\xi} k(I_j,s)\textup ds \biggr)\textup dx_v(\xi)-\overline{S}_j^v}\\
&\qquad \leq \varphi_1(1) \cdot 4\varepsilon^{-1} \sum_{j=1}^m \abs[\bigg]{(RS) \int_0^1 \biggl((L)\int_0^{\xi} k(I_j,s)\textup ds \biggr)\textup dx_v(\xi)-\overline{S}_j^v}\\
&\qquad \leq \frac{\varphi_1(1)}{1+\varphi_1(1)} \leq 1.
\end{align*} 
For the Riemann--Stieltjes sums we have
\begin{align*}
&\sum_{j=1}^m \varphi_j(4\varepsilon^{-1}\abs{\overline{S}_j^v - S_j^v})\\ 
&\qquad= \sum_{j=1}^m \varphi_j\Biggl(4\varepsilon^{-1}\abs[\Bigg]{\sum_{i=1}^n \sum_{l=\alpha_{i-1}+1}^{\alpha_i} \biggl( (L)\int^{\xi_i}_{\overline{\xi}_l} k(I_j,s)\textup ds\biggr)\bigl[x_v(\overline{\xi}_l)-x_v(\overline{\xi}_{l-1})\bigr]}\Biggr)\\
&\qquad\leq \sum_{j=1}^m \varphi_j\Biggl(4\varepsilon^{-1} \sum_{i=1}^n \sum_{l=\alpha_{i-1}+1}^{\alpha_i} \abs[\Bigg]{(L)\int^{\xi_i}_{\overline{\xi}_l} k(I_j,s)\textup ds}\abs[\big]{x_v(\overline{\xi}_l)-x_v(\overline{\xi}_{l-1})}\Biggr).
\end{align*}
As
\[
  \sum_{i=1}^n \sum_{l=\alpha_{i-1}+1}^{\alpha_i}\abs{x_v(\overline{\xi}_l)-x_v(\overline{\xi}_{l-1})} \leq \var x_v \leq 1,
\]
by the convexity and monotonicity of the Young functions we get
\begin{align*}
&\sum_{j=1}^m \varphi_j(4\varepsilon^{-1}\abs{\overline{S}_j^v - S_j^v})\\ 
&\qquad \leq \sum_{i=1}^n \sum_{l=\alpha_{i-1}+1}^{\alpha_i} \sum_{j=1}^m \varphi_j\biggl(4\varepsilon^{-1}\abs[\bigg]{(L)\int^{\xi_i}_{\overline{\xi}_l} k(I_j,s)\textup ds}\biggr) \cdot \abs{x_v(\overline{\xi}_l)-x_v(\overline{\xi}_{l-1})}\\
&\qquad\leq \sum_{i=1}^n \sum_{l=\alpha_{i-1}+1}^{\alpha_i} \var_\Phi \biggl(4\varepsilon^{-1} (L)\int^{\xi_i}_{\overline{\xi}_l} k(\cdot,s)\textup ds \biggr) \cdot \abs{x_v(\overline{\xi}_l)-x_v(\overline{\xi}_{l-1})}\\
&\qquad\leq \sum_{i=1}^n \sum_{l=\alpha_{i-1}+1}^{\alpha_i}\abs{x_v(\overline{\xi}_l)-x_v(\overline{\xi}_{l-1})} \leq 1.
\end{align*} 
Similarly, 
\begin{align*}
&\sum_{j=1}^m \varphi_j(4\varepsilon^{-1}\abs{S_j^v})\\
&\qquad=\sum_{j=1}^m \varphi_j\Biggl(4\varepsilon^{-1}\abs[\Bigg]{\sum_{i=1}^n \biggl( (L)\int_0^{\xi_i} k(I_j,s)\textup ds\biggr)\bigl[x_v(\xi_i)-x_v(\xi_{i-1})\bigr]}\Biggr)\\
&\qquad \leq \sum_{j=1}^m \varphi_j\Biggl(\sum_{i=1}^n\abs[\Bigg]{ \mu \cdot (L)\int_0^{\xi_i} k(I_j,s)\textup ds} \cdot 4\varepsilon^{-1}\mu^{-1}\abs[\big]{x_v(\xi_i)-x_v(\xi_{i-1})}\Biggr)\\
 &\qquad\leq\sum_{i=1}^n \sum_{j=1}^m \varphi_j\biggl(\abs[\bigg]{\mu \cdot (L)\int_0^{\xi_i} k(I_j,s)\textup ds}\biggr) \cdot 4\varepsilon^{-1}\mu^{-1}\abs[\big]{x_v(\xi_i)-x_v(\xi_{i-1})}\\
&\qquad\leq\sum_{i=1}^n \var_\Phi\biggl(\mu\cdot (L)\int_0^{\xi_i} k(\cdot,s)\textup ds\biggr) \cdot 4\varepsilon^{-1}\mu^{-1}\abs[\big]{x_v(\xi_i)-x_v(\xi_{i-1})}\\
 &\qquad\leq 4\varepsilon^{-1}\mu^{-1}\sum_{i=1}^n\abs[\big]{x_v(\xi_i)-x_v(\xi_{i-1})} \leq 1.
\end{align*}
Summing up all the partial estimates, we finally get $\sum_{j=1}^m \varphi_j\bigl(\varepsilon^{-1} \abs{(Kx_v)(I_j)}\bigr) \leq 1$. As the collection $I_1,\ldots,I_m$ is arbitrary, this implies that $\var_\Phi(\varepsilon^{-1} Kx_v) \leq 1$. Thus, $\abs{ Kx_v}_\Phi \leq \varepsilon$ for $v \geq v_0$, or equivalently, $\lim_{v \to \infty} \abs{Kx_v}_\Phi=0$. 

Using the Lebesgue dominated convergence theorem, we also have
\[
 \lim_{v \to \infty} (Kx_v)(0) = \lim_{v \to \infty} (L)\int_0^1 k(0,s)x_v(s)\textup ds=0.
\]
Therefore, $\lim_{v \to \infty} \norm{Kx_v}_\Phi=0$. To end the proof it now suffices to apply Proposition~\ref{prop:compactness_equiv}.
\end{proof}

\begin{remark}
Theorem~\ref{thm:K_BV_PBV_compact} extends an analogous result for linear integral operators acting in $BV[0,1]$ established in~\cite{PK}. We would also like to draw the readers' attention to the last section of that paper, where the author provided an in-depth discussion of various conditions guaranteeing compactness of linear integral operators acting into $BV[0,1]$ that can be found throughout the literature.   
\end{remark}

\begin{remark}
It turns out that conditions~\ref{H2} and~\ref{H3} have an interesting interpretation. Namely, let us consider the mapping $F \colon [0,1] \to \Phi BV[0,1]$ which to each point $\xi \in [0,1]$ assigns the function $t \mapsto (L)\int_0^\xi k(t,s)\textup ds$ of Schramm bounded variation. Then, under the general hypothesis~\ref{H1}, conditions~\ref{H2} and~\ref{H3} correspond, respectively, to boundedness and continuity of the mapping $F$. (For some comments concerning this interpretation see also~\cite{PK}*{Remark~3 and p.~22} and cf.~\cite{G}*{Section II.9}.)
\end{remark}

\end{document}